 \documentstyle[12pt]{article}

 \input amssym.def
\input amssym.tex
\begin{document}

\def \wt{{\rm wt}}
\def \tr{{\rm tr}}
\def \span{{\rm span}}
\def \Res{{\rm Res}}
\def \Res{{\rm QRes}}
\def \End{{\rm End}}
\def \E{{\rm End}}
\def \Ind {{\rm Ind}}
\def \Irr {{\rm Irr}}
\def \Aut{{\rm Aut}}
\def \Hom{{\rm Hom}}
\def \mod{{\rm mod}}
\def \ann{{\rm Ann}}
\def \<{\langle}
\def \>{\rangle}
\def \t{\tau }
\def \a{\alpha }
\def \e{\epsilon }
\def \l{\lambda }
\def \L{\Lambda }
\def \b{\beta }
\def \t{\theta}
\def \om{\omega }
\def \o{\omega }
\def \c{\chi}
\def \ch{\chi}
\def \cg{\chi_g}
\def \ag{\alpha_g}
\def \ah{\alpha_h}
\def \ph{\psi_h}
\def \be{\begin{equation}}
\def \ee{\end{equation}}
\def \bex{\begin{example}\label}
\def \bea{\begin{eqnarray} }
\def \eea{\end{eqnarray} }
\def \eex{\end{example}}
\def \bl{\begin{lem}\label}
\def \el{\end{lem}}
\def \bt{\begin{thm}\label}
\def \et{\end{thm}}
\def \bp{\begin{prop}\label}
\def \ep{\end{prop}}
\def \br{\begin{rem}\label}
\def \er{\end{rem}}
\def \bc{\begin{coro}\label}
\def \ec{\end{coro}}
\def \bd{\begin{de}\label}
\def \ed{\end{de}}
\def \pf{{\bf Proof. }}
\def \voa{{vertex operator algebra}}

\newtheorem{theorem}{Theorem}[section]
\newtheorem{corollary}{Corollary}[section]
\newtheorem{lemma}{Lemma}[section]
\newtheorem{conjecture}{Conjecture}[section]
\newtheorem{convention}{Convention}[section]
\newtheorem{remark}{Remark}[section]
\newtheorem{definition}{Definition}[section]
\newtheorem{proposition}{Proposition}[section]
\renewcommand{\baselinestretch}{1.4}
\newcommand{\binom}{ { { {\rm {wt}}  a} \choose i } }
\newcommand{\wta}{{\rm {wt} }  a }
\newcommand{\R}{\mbox {Res} _{z} }
\newcommand{\la}{\langle}
\newcommand{\ra}{\rangle}
\newcommand{\wtb}{{\rm {wt} }  b }
\newcommand{\g}{\bf g}
\newcommand{\hg}{\hat {\bf g} }
\newcommand{\h}{\bf h}
\newcommand{\hh}{\hat {\bf h} }
\newcommand{\nn}{\bf n}
\newcommand{\Z}{\Bbb Z}
\newcommand{\N}{\Bbb N}
\newcommand{\C}{\Bbb C}
\newcommand{\Q}{\Bbb Q}
\newcommand{\1}{\bf 1}
\newcommand{\NS}{\bf NS}
\newcommand{\cA}{\cal A}
\newcommand{\hn} {\hat {\bf n}}
\newcommand{\hf}{\frac{1}{2} }
\newcommand{\non}{\nonumber}
\def \l {\lambda}
\baselineskip=12pt
\newcommand{\rw}{\rightarrow}
\newcommand{\n}{\:^{\circ}_{\circ}\:}
\makeatletter
\@addtoreset{equation}{section}
\def\theequation{\thesection.\arabic{equation}}
\makeatother
\makeatletter

\title {Rationality of  unitary $N=2$ vertex operator superalgebras }

\author{ Dra\v zen Adamovi\' c }
\date{}
\pagestyle{myheadings}
\markboth
{Dra\v{z}en Adamovi\'{c}}
{}
\maketitle

\begin{abstract}
In this note we prove that the vertex operator superalgebras
associated to the unitary representations of the $N=2$
superconformal algebra are rational.
\end{abstract}

\section*{Introduction}

In this note we continue to study vertex operator superalgebras
(SVOAs) associated to the representations of $N=2$ superconformal
algebras. In \cite{A1} we presented the results which gave the
complete classification of irreducible modules. We proved that
SVOAs associated to the unitary $N=2$--representations have
finitely many irreducible modules. In non-unitary case the
representation theory is more complicated (and perhaps, more
interesting). The irreducible representations in non-unitary case
have a nice description as a union of finitely many rational
curves. We should mention that the results from \cite{A1}  were
motivated by analysis of Eholzer and Gaberdiel from \cite{EG}, and
from the series of papers \cite{FST}, \cite{ST}, \cite{SS}.

In the present paper we will make the first step in describing the
complete category of modules for these  SVOAs. Actually, we will
prove that the category of finitely-generated modules for SVOAs
associated to the unitary $N=2$--representations is semisimple.
This will imply that the corresponding SVOAs are rational. Our
proof of rationality is based on the analysis of embedding
diagrams of $N=2$ Verma modules from \cite{D} (see also
\cite{SS}), and the classification result from \cite{A1}.

 We also
believe that the main result of this note will be useful for the
problem of the classification of rational SVOAs with a unitary
central charges.

\section{$N=2$ superconformal vertex algebra}

In this section we recall the results from \cite{A1} on
representation theory of SVOAs associated to the $N=2$
superconformal algebra. We should also mention that
 the study of these SVOAs was initiated \cite{EG}.

$N=2$ superconformal algebra ${\cA}$ is the infinite-dimensional
Lie super algebra  with basis $ L_n, T_n, G_r ^{\pm}, C$, $n\in
{\Z}$, $ r\in {\hf} + {\Z}$ and (anti)commutation relations given
by \bea && [L_m,L_n] = (m-n) L_{m+n} + \frac{C} {12} (m ^ 3 -m)
\delta_{m+n,0} \nonumber \\ && [L_m, G_r ^{\pm}] =( {\hf } m -r )
G_{m+r} ^{\pm}  \non \\ &&[L_m, T_n] = - n T_{n+m} \non \\ &&
[T_m, T_n ] =  \frac {C}{3} m \delta _{m+n,0} \non \\ &&  [T_m,
G_r ^{\pm} ] = \pm G_{m+r} ^{\pm} \non \\ && \{ G_r ^{+}, G_s ^{-
} \} = 2 L_{r+s} + (r-s) T_{r+s} +
    \frac{C}{3} ( r  ^ 2 - \frac{1}{4} ) \delta_{r+s,0} \non \\
&& [L_m,C]= [T_n,C] = [G_r ^{\pm}, C] =0 \non \\
&& \{ G_r ^{+}, G_s ^{+} \}= \{ G_r ^{-}, G_s ^{-} \} = 0 \non
\eea
for all $m,n \in  {\Z}$, $r,s \in   {\hf} + {\Z}$.

We denote the Verma module generated from a highest weight vector
$\vert h,q,c \rangle$
with $L_0$ eigenvalue $h$, $T_0$ eigenvalue $q$ and central charge $c$
by $ M _{h,q,c}$. An element $v \in M _{h,q,c}$ is called singular
vector if
$$
L_n v = T_n v = G_r ^{\pm} v = 0,
\ \ n, r+ {\hf} \in {\N},
$$
and $v$ is an eigenvector of $L_0$ and $T_0$.
Let $J _{h,q,c}$ be the maximal $U({\cal A})$--submodule in  $ M _{h,q,c}$.
Then
$$L_{h,q,c} =\frac{M _{h,q,c}}{ J _{h,q,c} }$$
is the irreducible highest weight module.

Now we will consider the Verma module $M_{0,0,c}$.
One easily sees that  for every $c \in {\C}$
$$G_{ {-\hf} } ^{\pm} \vert 0,0,c \rangle $$
are the singular vectors in $M_{0,0,c}$.
Set
$$
V_c = \frac {M_{0,0,c} }
{ U( {\cA} ) G_{ {-\hf} } ^{+} \vert 0,0,c \rangle   +
 U( {\cA} ) G_{ {-\hf} } ^{-} \vert 0,0,c \rangle        } .
  $$
Then $V_c$ is a highest weight  ${\cA}$--module. Let ${\1}$ denote the
highest weight vector. Let $L_c=L_{0,0,c}$ be the corresponding simple module.
Define the following four vectors in $V_c$:
$$
\tau ^{\pm} = G_{-\frac {3}{2} } ^{\pm} {\1}, \
 j = T_{-1} {\1} , \
  \nu = L_{-2} {\1},
$$ and set \bea && G^{+} (z) = Y(\tau ^{+} ,z) = \sum _{n \in {\Z}
} G_{n+ {\hf} } ^{+} z ^{-n-2}, \non \\ && G^{-} (z) = Y(\tau ^{-}
,z) = \sum _{n \in {\Z} } G_{n+ {\hf} } ^{-} z ^{-n-2}, \non \\ &&
L (z) = Y(\nu,z) = \sum _{n \in {\Z} } L_{n }  z ^{-n-2}, \non \\
&& T (z) = Y(j,z) = \sum _{n \in {\Z} } T_{n }  z ^{-n-1}.
{\label{polja}} \eea It is easy to see that the fields  $G^{+}
(z)$, $G^{-} (z)$, $L(z)$, $T(z)$ are mutually local and the
theory of local fields (cf. \cite{K}, \cite{Li} ) implies the
following result.

\begin{proposition}
There is a unique extension of the fields (\ref{polja}) such that
$V_c$ becomes vertex operator superalgebra (SVOA).  Moreover,
$L_c$ is a simple SVOA.
\end{proposition}

We are interested in the list of all irreducible $L_c$--modules.
 It is very easy to see that every irreducible $L_c$--module has
to be irreducible highest weight $U({\cA})$--module. To verify
this statment, it is enough to notice that the Zhu' s algebra
$A(L_c)$ is a certain quotient of the polynomial algebra
${\C}[x,y]$ (cf. \cite{EG}, \cite{Z}).

 We will now present the classification result from \cite{A1}.

\begin{definition} A rational number $m=t/u$ is called admissible if
$u\in {\N}$, $t\in {\Z}$, $(t,u)=1$ and $2u+t-2 \ge 0$.
\end{definition}

For $m = \frac{t}{u}$ admissible set $c_m =\frac{3m}{m+2}$, $N= 2
u + t -2$ and $$ S^m = \{ n - k (m+2) \ \vert \ k,n \in {\Z}_+,
n\le N, k\le u-1\}. $$ Let $$ W^{ c_m} = \left \{ \left( \frac{j k
- \frac{1}{4} }{m+2},
 \frac{j-k}{m+2} \right) \ \vert \ j,k \in {\N} _{ \frac{1}{2} },
  0< j,k, j+k \le N+1 \right \},
$$ and if $m \notin {\N}$ let $$ D^{c_m} = \{ (h,q) \in {\C} ^2 \
\vert \ q^2 + \frac{4h}{m+2} = \frac{r (r+2)}{ (m+2) ^2}, \ r\in
S^m \setminus {\Z} \}. $$ Note that the set $W ^{ c_m}$ is finite
and the set $D^{c_m}$ is union of finitely many rational curves.

\begin{theorem} \cite{A1} \label{clas.unit}
Assume that $m \in {\N}$ and $c = c_m$.  Then the set $$
\{L_{h,q,c} \ \vert \  (h,q) \in W^ {c_m} \} $$ provides all $L_c$
irreducible modules for the SVOA $L_c$. So, irreducible
$L_c$--modules are exactly all unitary modules for $N=2$
superconformal algebra with the central charge $c$.
\end{theorem}

\begin{remark}
Theorem \ref{clas.unit} shows that SVOA $L_{c_m}$ for $m \in {\N}$
has exactly $\frac{(m+2) (m+1)} {2} $ non-isomorphic irreducible
modules.
\end{remark}
\begin{theorem} \cite{A1} \label{clas.nonunit}
Assume that $m \in {\Q}$ is admissible such that $ m \notin {\N}$. Let
$c = c_m$. Then the
set
$$
\{ L_{h,q,c} \ \vert \  (h,q) \in  W^ {c_m} \cup D^ {c_m} \}
$$
provides all $L_c$ irreducible modules for the
SVOA $L_c$.
\end{theorem}

 \section{Rationality of the SVOA $L_{ c_m}$ for $m \in {\N}$ }

In this section we will present a proof of the rationality of the
SVOA $L_{ c_m}$. Theorem \ref{clas.unit} gives that SVOA $L_{
c_m}$ has finitely many irreducible modules. We will now finish the proof
of rationality by proving that every finitely generated $L_{c_m}$--module
is completely reducible.
 In almost all known cases the proof of 
rationality is based on the structure theory of certain categories
of highest weight representations, for example the category $\cal
O$ in the case of affine Lie algebras (cf. \cite{DGK}). In order
to modify the proof from affine case to superconformal algebras,
one has to know the embedding diagrams for Verma modules. In the
case of unitary representations of $N=2$ superconformal algebra,
the correct embedding diagrams for Verma modules was done by M.
Dorrz\"apf in \cite{D}. He also proved that unitary $N=2$ Verma
modules do not contain subsingular vectors, which implies that
every submodule of a unitary Verma module is generated by singular
vectors (cf. \cite{D} and the corresponding references).

\begin{definition}  SVOA $V$ is called {\bf rational} if $V$ has finitely
many irreducible modules, and if every finitely generated
$V$--module is completely reducible.
\end{definition}

In this section, let $m \in {\N}$, and $c=c_m = \frac{3m}{m+2}$.
The results from \cite{D}  imply the following analysis.

 For $(h,q) \in W^{c}$, we have that
$$h = \frac{j k - \frac{1}{4} }{m+2}, q= \frac{j-k}{m+2}, $$
 for
$$j,k \in {\N} _{ \frac{1}{2} }, 0< j,k, j+k \le m+1. $$
We are interested in the weights of all homogeneous singular
vectors in Verma module $M_{h,q,c}$.

  Theorem 4.A from \cite{D} implies that the Verma module $M_{h,q,c}$ is
reducible, and has infinitely many  singular vectors
\bea \label{listv}
 v_{h_1 ^{i}, q_1 ^{i}}, v_{h_2 ^{i}, q_2 ^{i}},
\dots , v_{h_6 ^{i}, q_6 ^{i}},  \quad i\in {\N} \eea
of the weights
\bea \label{tezine}
&& h_1 ^{i} = h + i ( i (m+2) -j -k ), \quad q^{i} _1 =q,
\label{t1}
\\
&& h_2 ^{i} = h + i ( i (m+2) +j +k ), \quad q^{i} _2 =q,
\label{t2}
\\
&& h_3 ^{i} = h + k + (i-1) ( i (m+2) +j +k ), \quad q^{i} _3 =q
+1, \label{t3} \\
 && h_4 ^{i} = h + k+ (i+1) ( i (m+2) -j -k ),
\quad q^{i} _4=q+1, \label{t4} \\
&& h_5 ^{i} = h + j + (i-1) ( i (m+2) +j +k ), \quad q^{i} _5 =q
-1, \label{t5} \\
&& h_6 ^{i} = h + j + (i+1) ( i (m+2) -j -k ), \quad q^{i} _6 =q
-1. \label{t6} \eea
Moreover, the maximal submodule $M^{1} _{h,q,c}$ of the Verma
module $M_{h,q,c}$ is generated by three vectors :

$$ v_{h_1 ^{1}, q_1 ^{1}},  v_{h_3 ^{1}, q_3 ^{1}},  v_{h_5 ^{1},
q_5 ^{1}}. $$

We have the following lemma.

\begin{lemma} \label{duv}
 Assume that $v_{h',q'}$ is a singular vector of weight $(h',q')$ in the
Verma module $M_{h,q,c}$, where $c=c_m$, $(h,q) \in W^{c}$. Then
$(h',q') \notin W^{c}$.
\end{lemma}
{\em Proof.} Since the weights (\ref{t1})-(\ref{t6})  provide the
weights of all singular vectors in the Verma module $M_{h,q,c}$,
we have to check the statement of the Lemma for every singular
vector in the list (\ref{listv}). For simplicity,  we will only
consider vectors $ v_{h_1 ^{i}, q_1 ^{i}}$. For other singular
vectors, the considerations are completely analogous. So, we will
prove that
for every $i \in {\N}$
 $$ (h_1 ^{i}, q_1 ^{i}) \notin W^{c}.$$
 First, we notice that
 $$h_1 ^{i} = \frac{(j - i (m+2) ) ( k- i (m+2) ) -
 \frac{1}{4}}{m+2}, \quad q_1 ^{i} = \frac{j-k}{m+2}.$$
 Assume that $({h_1 ^{i}, q_1 ^{i}}) \in W^{c}$. Then there are
\bea \label{uvjet}
 \overline{j},\overline{k} \in {\N} _{ \frac{1}{2} }, 0< \overline{j},\overline{k},
 \overline{j} +\overline{k} \le m+1,
 \eea
 such that
 \bea \label{sustav}
 && h_1 ^{i} = \frac{ \overline{j} \  \overline{k} -\frac{1}{4}
 }{m+2}, \quad q_1 ^{i} = \frac{\overline{j} - \overline{k}}{m+2}.
 \eea
 But, the  equation (\ref{sustav}) implies that
 $$\overline{j} = j -i (m+2), \overline{k} = k-i (m+2)$$
 or
$$\overline{j} = -j +i (m+2), \overline{k} = -k+i (m+2),$$ which
in both cases contradicts the condition (\ref{uvjet}).${ \ \
\Box}$ \vskip 10mm

 Lemma \ref{duv} and the fact that unitary Verma
modules don't contain  subsingular vectors imply the following
proposition.

\begin{proposition} \label{lnotin}
Let $V_{h,q,c}$ be any highest weight $U({\cal A})$--module such
that the highest weight $(h,q) \in W^c$, $c = c_m$. Assume that
$V_{h,q,c}$  is reducible. Then there is a highest weight $
U({\cal A})$--submodule $N'$ with the highest weight $(h',q')
\notin W^c$.
 \end{proposition}

Proposition \ref{lnotin} and Theorem \ref{clas.unit} imply the
following.
 \begin{corollary} \label{lnotin2}
 Assume that $M$ is a highest weight  $U({\cal A})$--module  such that
 $M$ is a module for SVOA $L_{ c_m}$, $m \in {\N}$. Then $M$ is irreducible.
 \end{corollary}

\begin{remark} Proposition \ref{lnotin} and Corollary \ref{lnotin2} can
be also proved by using the results of A. M. Semikhatov and V. A.
Sirota from \cite{SS}. Their approach  uses the relations between
the embedding diagrams of $N=2$--modules and modules for affine
Lie algebra $\hat {sl_2}$, and it can be applied on non-unitary
highest weight $N=2$--modules. We hope to study the non-unitary
case in our next  publications (cf. \cite{A3}).
\end{remark}

 Let $\omega$ be the involutory
anti-automorphism of ${\cA}$, defined with $c \mapsto c$ and

 \bea
 &&L_m \mapsto L_{-m}, \quad T_m\mapsto T_{-m}, \nonumber \\
 &&  G^{+}_{ m-{\hf}} \mapsto  G^{-}_{- m+{\hf}}, \quad
G^{-}_{ m-{\hf}} \mapsto  G^{+}_{- m+{\hf}} \nonumber \eea
  for every $m \in {\Z}$.

 We will now prove that every finitely generated $L_c$--module is
 completely reducible. It is enough to consider $L_c$--modules
 with the following graduation:
 \bea
 &&M= \oplus _{ m
\in {\hf} {\Z}_+ } M_{ m}, \quad\mbox{dim} \ M_{m} < \infty
\nonumber \\
&& L_0 \ \mbox{ acts semisimply on} \ \  M_{m}. \nonumber \eea

 We define \cite{DGK}
$M^{{\omega}}=\oplus_{ m \in {\hf} {\Z}_+}M_{m}^{*}$ with the
following action $(af)(u)=f({\omega}(a)u)$ for any $f\in
M^{{\omega}}, a\in {\cA}, u\in M$. In the same way as in
\cite{DGK} one can see that  $M \mapsto M^{\omega}$ is a
contravariant functor and $$(M^{{\omega}})^{\omega}\simeq M,\;
L_{h,q,c}^{{\omega}}\simeq L_{h,q,c} \;\;\;\mbox{ for any } (h,q)
\in {\C} ^{2}.$$

 \begin{theorem}  \label{rat.unit}
 For $m\in {\N}$ the SVOA $L_{ c_m}$ is rational.
 \end{theorem}
 We proved in Theorem \ref{clas.unit}  that SVOA $L_{ c_m}$ has finitely many irreducible
 modules. For a proof of rationality, it remains  to prove that every finitely generated
 $L_{ c_m}$--module is completely reducible.  To see this, it is
 enough to prove the following lemma (see also \cite{DLM},
 \cite{A2}, \cite{W} for  proofs of rationality in some  similar
 cases).

 \begin{lemma} \label{ext} Let $c=c_m$ for $m \in {\N}$, and let
 $(h_1,q_1), (h_2,q_2) \in W^{c}$.
 Then every short exact sequence of $U({\cA})$--modules \bea
\label{seq} 0 \longrightarrow L_{h_1,q_1,c}
    \stackrel{\iota}{\longrightarrow} M
    \stackrel{\pi}{\longrightarrow} L_{h_2,q_2,c}
    \longrightarrow 0
\eea splits.
\end{lemma}
{\em Proof.} Without loss of generality, we may assume that $ h_2
\le h_1 $. Otherwise we can apply the  functor $ M\mapsto
M^{\omega}$ to the short exact sequence to reverse it. Let
$v_{h_2,q_2}$ be the highest weight vector of $L_{h_2,q_2,c}$.
Pick a vector $v_{h_2,q_2 }' \in M $ of weight $(h_2,q_2)$ such
that $\pi ( v_{h_2,q_2}') = v_{h_2,q_2}$. We claim that
$v_{h_2,q_2}'$ is a singular vector in $M$. Let $n \in {\N}$.
 Assume that
$T_{ n} . v_{h_2,q_2}' \neq 0 $ . Then $$\pi (T_{n} v_{h_2,q_2}')
= T_{n} \pi (v_{h_2,q_2}') = T_{n} v_{h_2,q_2} = 0.$$ This implies
\begin{equation}
T_{n} v_{h_2,q_2}' = \iota (u) \label{jed1}
\end{equation}
for some nonzero $u \in L_{h_1,q_1,c}$, since the short sequence
is exact. Comparing the weights of both sides of equation
(\ref{jed1}), we have $ h_2 - n =  h_1 + \alpha$ for nonzero
$\alpha \in {\hf} {\Z}_+$. It follows that $ h_2 = h_1 +\alpha + n
> h_1,$ which is a contradiction. In the same way the assumptions $
G^{\pm}_{n -{\hf} } v_{h_2,q_2}' \neq  0$, $L_{n}  v_{h_2,q_2}'
\neq  0$   for $n > 0$ lead to contradiction.

Denote by $M'$ the highest weight submodule of $M$ generated by
the singular vector $v_{h_2,q_2}'$. Assume that $M'$ is reducible.
Then Lemma \ref{lnotin} implies that $M'$ has a highest weight
submodule $N' \ne 0$ with the highest weight $(h',q') \notin
W^{c_m}$. Then irreducibility of $L_{h_2,q_2,c}$ implies that $N'
\subseteq \mbox{Ker} (\pi)$, which contradicts the fact that short
sequence is exact. In this way we have proved that $M'$ is
irreducible. This implies that $ M \cong M' \oplus L_{h_1,q_1,c}$,
and we have proved that the sequence (\ref{seq}) is exact. ${ \ \
\Box}$ \vskip 10mm

 Let us conclude this section by giving one
characterization of the category of $L_{ c_m}$--modules. We know
from the embedding structure that the maximal submodule of the
Verma module $M_{0,0,c}$ is generated three vectors. Two such
vectors are
 $G_{-\hf} ^{\pm} \vert 0,0,c \rangle$, and the third vector has
relatively charge zero and level $m +1$. In particular, we have
that the maximal submodule of $U(\cA)$--module $V_{ c_m}$ is
generated by one singular vector $v_{sing}$ of level $m+1$. Moreover,  the
results from \cite{ST}, \cite{FST} imply that the maximal
submodule of $V_{c_m}$ is generated by the vector
 $G^{-} _{-m- \frac{3}{2}} \cdots  G^{-} _{- \frac{3}{2}} {\1}$
 (or equiv.  $G^{+} _{-m- \frac{3}{2}} \cdots  G^{+} _{- \frac{3}{2}}
 {\1}$),
 and that the singular vector  is given by the following formulae:
$$v_{sing} = G^{+} _{ \frac{1}{2}} \cdots  G^{+} _{m+ \frac{1}{2}}
 G^{-} _{-m- \frac{3}{2}} \cdots  G^{-} _{- \frac{3}{2}} {\1}.
 $$

The previous analysis and the vertex operator superalgebra
structure on $V_{c_m}$ and $L_{ c_m}$ imply the following lemma.
\begin{lemma} \label{kar}
Let $M$ be a $U({\cA})$--module of central charge $c_m$. Then the
following statements are equivalent:
\item[(i)] $M$ is a $L_{c_m}$--module,
\item[(ii)] $Y(G^{-} _{-m- \frac{3}{2}} \cdots  G^{-} _{- \frac{3}{2}}
{\1},z) M = 0$, (or equiv. $Y(G^{+} _{-m- \frac{3}{2}} \cdots
G^{+} _{- \frac{3}{2}} {\1},z) M = 0$),
\item[(iii)] $ : \partial^{(m)} G^{-} (z) \cdots G^{-} (z) : M =
0$ \\(or equiv. $ : \partial^{(m)} G^{+} (z) \cdots G^{+} (z) : M
= 0)$.
\end{lemma}

Now Lemma \ref{kar} and Theorem \ref{rat.unit} imply the
following.

\begin{corollary} \label{kar2}
Assume that $M$ is a finitely generated $U({\cA})$--module of the
central charge  $c_m$ such that $ : \partial^{(m)} G^{-} (z)
\cdots G^{-} (z) : M = 0.$
 Then $M$ is a direct sum of unitary representations with the
 central charge $c_m$.
\end{corollary}

\begin{remark} The statment of the Corollary \ref{kar2} doesn't
need vertex algebra language, but our proof uses vertex algebras.
This statement was also known to physicist, and it was used in
\cite{FS} for some cohomology calculation.
\end{remark}
\begin{remark} Let $e,f,h$ be the generators of $sl_2$. It is
well-known that the integrable $\hat{sl_2}$--modules of level $m$
are characterized by the field-relation $e(z) ^{m+1}=0$, and the
Corollary \ref{kar2} gives $N=2$  interpretation of this fact. It
will be also  interesting to translate some other properties of
the representation theory of $\hat{sl_2}$ in $N=2$ language.
\end{remark}

Department of Mathematics, University of Zagreb,
Bijeni\v{c}ka 30, 10000 Zagreb, Croatia

E-mail address: adamovic@cromath.math.hr

\end{document}